\documentclass[12pt]{article}
\newcommand{\fr}{ \frac}
\newcommand{\beq}{\begin{equation}}
\newcommand{\eeq}{\end{equation}}
\newcommand{\br}{\begin{array}}
\newcommand{\er}{\end{array}}
\newcommand{\lb}{\label}
\newcommand{\ar}{\rightarrow}
\newcommand{\dl}{\delta}

\newcommand{\lan}{\langle}
\newcommand{\rn}{\rangle}
\newcommand{\pr}{\prime}

\newcommand{\lm}{\lambda}

\begin{document}
 \begin{center}
 {\Large\bf   Summation Formulas for the product of the q-Kummer Functions
 from  $E_q(2)$ }
 \end{center}
 \vspace{5mm}
 \begin{center}
 H. Ahmedov$^1$  and I. H. Duru$^{2,1}$
 \end{center}

 \noindent
 1.Feza G\"ursey Institute,  P.O. Box 6, 81220,
 \c{C}engelk\"{o}y, Istanbul, Turkey \footnote{E--mail :
  hagi@gursey.gov.tr and duru@gursey.gov.tr}.

 \noindent
 2.Trakya University, Mathematics Department, P.O. Box
 126, Edirne, Turkey.

 \vspace{5mm}
 \begin{center}
 {\bf Abstract}
 \end{center}

 \noindent
 Using  the representation of $E_q(2)$ on the non-commutative space
 $zz^*-qz^*z=\sigma$; $q<1$, $\sigma>0$   summation formulas for
 the product  of two, three and four q-Kummer functions are  derived.

 \vspace{1cm}
 \noindent
 {\bf 1. Introduction}

 \vspace{5mm}
 \noindent
 Properties of manifolds can be investigated by means of
 their automorphism groups. Non-commutative  spaces too are studied similarly.
 For example the quantum groups $E_q(2)$ and $SU_q(2)$
 are the symmetry groups of the quantum plane and the  quantum sphere
 respectively \cite{pod}-\cite{wor1}. The group representation theory
gives the
 possibility to construct  the complete set  of orthogonal functions on these deformed spaces.
 For example the Hahn-Exton  q-Bessel and q-Legendre functions
 appears as the matrix elements  of the unitary representations of $E_q(2)$
 \cite{wor}-\cite{koel}  and  $SU_q(2)$ \cite{koor} \cite{vil} which are
the complete
 set of  orthogonal functions  on the quantum plane  and the quantum sphere
 respectively.  Using group theoretical methods the invariant distance
 and the Green functions  have also been written on the quantum sphere
 \cite{duru1} and the quantum plane \cite{duru2}.

 In recent works we have studied the non-commutative space $[z,
 z^*]=\sigma$ ( i. e. the space  generated by the
 Heisenberg algebra) by means of its automorphism groups $E(2)$
 and $SU(1,1)$ \cite{ahm}, \cite{ahm1}. The basis in this noncommutative space
 where irreducible representations of $E(2)$ are realized were
 found to be the Kummer functions which involves the coordinates
 $z, \ z^*$ not as their arguments but as indices. That study
 enable us to obtain  generic summation formulas involving Kummer
 and Bessel functions. For $SU(1,1)$ case the basis is given in terms
 of the hypergeometric functions having the non-commutative coordinates
 $z$ and $z^*$ as the parameters. Again we derived generic summation formulas
 involving hypergeometric and Jacobi functions.
 This analysis enable us to construct different complete sets of orthogonal functions
 on the non-commutative space. Both studies also provide new group theoretical
 interpretations
 for the already known relations involving special functions.

 Motivated by the outcomes of the above mentioned studies,
 in the present work we consider the two parametric
 deformation of the plane which is the $*$-algebra $P^\sigma_q$
 generated by  $z$ and $z^*$ with
 \begin{equation}\lb{1com}
 zz^* -q z^* z = \sigma, \ \ \ \ q<1, \ \sigma > 0
 \end{equation}
 which possesses  the symmetry of the group $E_q(2)$.
 In $\sigma\ar 0$ limit it becomes the usual quantum plane.
 In $q\ar 1$ limit  it turns to  the algebra of functions on
 the Heisenberg algebra. This study allows us to obtain many
 identities involving several Hahn-Exton q-Bessel and Moak q-Laguerre
 functions which are the special forms of the q-Kummer functions.
 Note that previously some formulas involving q-Laguerre functions were derived by
 making use of the  representation  theory of the q-oscillator algebra
 \cite{ask1}-\cite{kaln3}. Some relations involving
 the basic Bessel and  Laguerre functions were also considered in
 \cite{koelnik}.

 In Section 2 we realize $E_q(2)$  as the automorphism group of
 the non-commutative space $P^\sigma_q$.

 In Section 3 we construct the basis in $P^\sigma_q$
 where the irreducible representations of $E_q(2)$ are realized.

 Section 4 is devoted to the generic summation formulas  for the product of two, three and four
 q-Kummer functions.

 In Section 5 some simple examples are presented.

 \vspace{1cm}
 \noindent
 {\bf 2. $E_q(2)$ as the symmetry group of $P^\sigma_q$}

 \vspace{2mm}
 \noindent
 Euclidean group transformations of $z, \ z^*$ plane are given by
 \begin{eqnarray}\lb{da}
 \delta(z) & = & B + A z, \\
 \delta(z^*) & = & B^* + A^* z^*.
 \end{eqnarray}
 The conditions
 \beq\lb{1oss1}
 \delta(z)\delta (z^*) -q \delta (z^*) \delta (z) = 1
 \eeq
 and
 \beq\lb{1st}
 (\delta(z))^* = \delta (z^*)
 \eeq
 imply the $E_q(2)$ relations
 \beq\lb{1c}
 BB^* = qB^*B, \ \  AB=qBA, \ \ AB^*=qB^*A, \ \ A^*=A^{-1}.
 \eeq
 Note that  $z, \ z^*$ commute with $B, \ B^*$ and $A$.
 (Starting
 with(\ref{1oss1}) we employ $\sigma =1$. When we need to
 calculate $\sigma \ar 0$ limit we replace $z, z^*$ by
 $\fr{z}{\sigma}, \fr{z^*}{\sigma}$). Formulas
 \begin{eqnarray}\lb{1oss}
 z|n, j\rn & = & \sqrt{(n)_q} |n-1, j\rn, \nonumber \\
 z|n, j\rn & = & \sqrt{(n+1)_q} |n+1, j\rn \nonumber \\
 B|n, j\rn & = & q^{\fr{j}{2}}|n, j-1\rn, \\
 B^*|n,_j\rn & = & q^{\fr{j+1}{2}}|n, j+1\rn, \nonumber \\
 A|n, j\rn & = &|n, j-2\rn,\nonumber
 \end{eqnarray}
 where
 $$
 (n)_q = \fr{1-q^n}{1-q}
 $$
 give the solution of (\ref{1com}) and (\ref{1c}) in some suitable domain of the
 Hilbert space $X$ with the basis $\{|n, j\rn\}$, $n=0, 1, 2,
 \dots$ and0, $j\in \bf{Z}$.
 Let us define in $X$ a new basis such that
 \begin{eqnarray}\lb{dar}
 \delta(z) |n, j\rn^\pr & = &\sqrt{(n)_q} |n-1, j\rn^\pr \\
 \delta(z) |n, j\rn^\pr & = &\sqrt{(n)_q} |n +1, j\rn^\pr.
 \end{eqnarray}
  Due to
 \begin{equation}
 z e_q^{-x z^*}= - x e_q^{x z^*}+e_q^{-q x z^*}z,
 \end{equation}
 with
 \begin{equation}
 e_q^{x}= \sum_{k=0}^\infty \fr{x^k}{(k)_q!}
 \end{equation}
 being the q-deformed exponential function, we conclude that
 \beq
 |0, j\rn^\pr = e_q^{-A^*B a^*} \sqrt{e_q^{-B^*B}} |0, j\rn
 \eeq
 is the ground state of the new basis:
 \beq
 \delta(z)|0, j\rn^\pr = 0.
 \eeq
 Applying the creation operator $(\delta (z))^*$ on this state we
 can generate the desired basis in $X$:
 \beq
 |n, j\rn^\pr = \fr{(\delta(z^*))^n}{\sqrt{(n)_q!}} |0, j\rn^\pr.
 \eeq
 We also have
 \beq\lb{uau}
 \delta(z)= U z U^*,
 \eeq
 where  $U$ is the unitary operator in
 \beq
 |n, j\rn^\pr =U |n, j\rn.
 \eeq
 Thus $\dl$ defines the homomorphic map of $P^\sigma_q$ into the
 $*$-algebra generated by $A, \ B, \ z$ and their adjoints.
 Applying twice this map we get
 \beq
 \delta^\pr(\delta(z))= B^\pr + A^\pr \delta(z)=B^{\pr\pr}
 +A^{\pr\pr} z,
 \eeq
 where
 \beq
 B^{\pr\pr} = B^\pr + A^\pr B, \ \ \ \
 A^{\pr\pr} = A^\pr A,\ \
 \eeq
 is the group multiplication in $E_q(2)$ ( the operator $B^{\pr\pr}$ has  the same analytic
 properties as $B$   \cite{wor2}). $\delta$  defines the representation of $E_q(2)$ in
 $P^\sigma_q$.

 Before closing this section we give the explicit formula for the
 matrix representation of $U$:
 \beq
 U_{(mi)(nj)} = \lan m i| n, j\rn^\pr.
 \eeq
 For $|n, j\rn =|n\rn|j\rn$ we first define
 \beq
 U_{mn} = \lan n |\delta(z^*))^n  e_q^{-A^*B z^*} |0\rn
 \sqrt{\fr{e_q^{-B^*B}}{(n)_q!}}
 \eeq
 which is the function of $B$, $B^*$, $A$ and  $A^*$.
 Then
 \beq
 U_{(mi)(nj)} =  \lan i | U_{mn} |j\rn.
 \eeq
 After some algebra  we get
 \beq\lb{matru1}
 U_{mn} = A^{-m}B^{* n-m}\Phi_{mn}(\eta) \ \ \ \mathrm{for} \ n\geq m
 \eeq
 and
 \beq\lb{matru2}
 U_{mn} = q^{\fr{(m-n)(m-n-1)}{2}}
 A^{-m}(-B)^{m-n}\Phi_{nm}(\eta) \ \ \ \mathrm{for} \ m\geq n,
 \eeq
 where $\eta^2\equiv B^*B$ and
 \beq
 \Phi_{mn}(\eta) = \sqrt{\fr{(n)_q!}{(m)_q!}}
 \fr{\sqrt{e_q^{-\eta^2}} }{(n-m)_q!}
 \Phi^q (q^{-m}, q^{1+n-m}; q^{n+1}\eta^2).
 \eeq
 Here
 \beq
 \Phi^q (a, b; x)=
 \sum_{k=0}^\infty \fr{q^{\fr{k(k-1)}{2}}(a; q)_k}{(q; q)_k (b; q)_k}
 ((1-q)x)^k
 \eeq
 which in $q\ar 1$ limit  reduces to the  Kummer function:
 \beq
 \lim_{q\ar 1} \Phi^q (q^c, q^d; x) = \Phi (c, d; x).
 \eeq
 We call it the q-Kummer function. The functions $\Phi_{nm}$ can
 also be expressed in terms of  Moak's q-Laguerre polynomials
 \cite{moak}
 \beq
 L^{q(\alpha)}_n(x) = \fr{(q^{1+\alpha}; q)_n}{(q; q)_n}
 \Phi^q (q^{-n}, q^{1+\alpha}; q^{1+\alpha +n}x)
 \eeq
 as
 \beq
 \Phi_{mn}(\eta) = \sqrt{e_q^{-\eta^2}\fr{(m)_q!}{(n)_q!}}
 L^{q(n-m)}_m (\eta^2), \ \ \ \mathrm{for} \ n\geq m.
 \eeq

 \vspace{1cm}
 \noindent
 {\bf 3. Irreducible representations of $E_q(2)$ in $P^\sigma_q$}

 \vspace{2mm}
 \noindent
 The irreducible representation of the deformed enveloping algebra
 $U_q(e(2))$
 \beq
 P^*P=q P P^*, \ \ \ KP=qPK, \ \ \ P^*K=q KP^*
 \eeq
 defined by the weight $\lm\in \bf{R}$ in the non-commutative space
$P_q^\sigma$
 is given by
 \beq\lb{3rep}
 R(K) D^\lm_j(z,z^*)  =  q^j  D^\lm_j(z, z^*)
 \eeq
 \beq\lb{3repa}
 R(P) D^\lm_j(z,z^*)  =  \lm q^{\fr{j}{2}}D^\lm_{j-1}(z,z^*),
 \eeq
 \beq\lb{3repb}
 R(P^*) D^\lm_j(z,z^*) =  \lm q^{\fr{j+1}{2}}D^\lm_{j+1}(z,z^*),
 \eeq
 where
 \begin{eqnarray}
 R(K) z^{*n}z^m & = & q^{m-n}z^{*n}z^m \nonumber \\
 R(P) z^{*n}z^m & = & iq^{-n}(m)_q z^{*n}z^{m-1} \\
 R(P^*) z^{*n}z^m & = & iq^{-n+1}(n)_q z^{*n-1}z^m \nonumber
 \end{eqnarray}
 defines the right realization of $U_q(e(2))$.
 (\ref{3rep}) implies
  \beq\lb{ans}
  D^\lm_j(z, z^*) =\{
  \begin{array}{c}
   f^\lm_j(\zeta )z^j \ \ \ \mathrm{for }\ j\geq 0 \\
   z^{*-j} f^\lm_{-j}(\zeta ), \ \ \ \mathrm{for} \ j\leq 0
  \end{array}
  \eeq
  where
  \beq
  \zeta \equiv 1-(1-q)z^*z.
  \eeq
 Inserting the ansatz (\ref{ans}) in (\ref{3repa}) and (\ref{3repb})
 we get
 \beq\lb{for}
 f^\lm_j(\zeta ) = \fr{ q^{\fr{j^2}{4}} (i\lm)^j} {(j)_q!}
 \Phi^q (\zeta^{-1}, q^{j+1}; \ q^{j+1}\lm^2\zeta ).
 \eeq
 In the derivation of (\ref{for}) we used
 \beq
 z^{*n}z^n = (-)^n q^{\fr{n(1-n)}{2}}(\zeta^{-1}; q)_n \zeta^n.
 \eeq
 By means of the universal T-matrix we can exponentiate (\ref{3rep})
 and get
 \beq\lb{add}
 \delta (D^\lm_j) = \sum_{i=-\infty}^\infty t^\lm_{ji} D^\lm_i,
 \eeq
 where
 \begin{eqnarray}\lb{matr}
 t^\lm_{ij} & = & q^{\fr{i^2-j^2}{4}}\fr{
(i\lm B)^{i-j}A^j}{(i-j)_q!}
 \Phi^q ( 0, q^{1+i-j}; (q-1) q^{1-j} (\lm\eta)^2) \ \
 \mathrm{for} \ i\geq j \\
 t^\lm_{ij} & = & q^{\fr{i^2-j^2}{4}}\fr{
 A^j(i\lm B^*)^{j-i}}{(j-i)_q!}
 \Phi^q ( 0, q^{1+j-i}; (q-1) q^{1-i} (\lm\eta)^2) \ \ \mathrm{for} \ j\geq i
 \end{eqnarray}
 are the matrix elements of the irreducible  representations of
 $E_q(2)$ \cite{bon}, \cite{koel}. We can express them in terms of
 the Hahn-Exton q-Bessel functions \cite{sw}
 \beq
 J_k^q(x) =\fr{x^k}{(k)_q!}
 \Phi^q ( 0; q^{1+k} |q; (q-1)qx^2)
 \eeq
 as
 \begin{eqnarray}
 t^\lm_{ij} & = &
 (\sqrt{-1}q^{\fr{1}{4}})^{i-j}
 V^{i-j}A^j J_{i-j}^q( q^{-\fr{j}{2}}\lm\eta)\ \ \ \mathrm{for} \
 i\geq j \\
 t^\lm_{ij} & = &
 (\sqrt{-1}q^{-\fr{1}{4}})^{j-i}
 V^{i-j}A^j J^q_{j-i}( q^{-\fr{i}{2}}\lm\eta)\ \ \ \mathrm{for} \  j\geq i,
 \end{eqnarray}
 where $V$ is the unitary operator defined by $B=V\eta$.

 In $\sigma \ar 0$ limit the non-commutative space $P_q^\sigma$
 becomes the quantum plane $E_q(2)/U(1)$ generated by $B$, $B^*$:
 \beq
 \lim_{\sigma\ar 0}D_j^{\sqrt{\sigma}\lm}
 (\fr{B}{\sqrt{\sigma}}, \fr{B^*}{\sqrt{\sigma}}) =
 t^\lm_{j0}.
 \eeq

 In $q \ar 1$ limit $P_q^\sigma$ becomes the noncommutative space
 generated by the Heisenberg algebra  \cite{ahm}:
 \beq
 \lim_{q\ar 1}D_j^{\lm}(z, z^*) =
 \Phi(-zz^*; \ 1+j; \lm^2)\fr{(i\lm z)^j}{j!}.
 \eeq

 In $\sigma \ar 0$ and $q\ar 1$  limit  we arrive at  the  complex plane $E(2)/U(1)$:
 \beq
 \lim_{\sigma\ar 0}\lim_{q\ar 1}D_j^{\sqrt{\sigma}\lm}
 (\fr{re^{i\psi}}{\sqrt{\sigma}}, \fr{re^{-i\psi}}{\sqrt{\sigma}})
 = i^je^{-ij\psi} J_j(\lm r).
 \eeq

 \vspace{1cm}
 \noindent
 {\bf 4. Summation formulas for the q-Kummer functions}

 \vspace{2mm}
 \noindent
 (\ref{uau}) and (\ref{add}) imply
 \beq\lb{41}
 UD^\lm_jU^* = \sum_{i=-\infty}^\infty t^\lm_{ji} D^\lm_i
 \eeq
 or
 \beq\lb{42}
 UD^\lm_j = \sum_{i=-\infty}^\infty t^\lm_{ji} D^\lm_iU,
 \eeq
 \beq\lb{43}
 D^\lm_j = \sum_{i=-\infty}^\infty t^\lm_{ji} U^*D^\lm_i U.
 \eeq
 The above formulas define the summation of products of two, three
 and four q-Kummer functions.
 Sandwiching (\ref{41}), (\ref{42}) and (\ref{43}) between the states $\lan m|$ and
 $|n\rn$ and using
 \beq
 (D^\lm_j)_{mn}=\sqrt{\fr{(n)_q!}{(m)_q!}}f^\lm_j(q^m)
 \dl_{j,n-m}, \ \ \  \mathrm{for} \  j\geq 0,
 \eeq
 \beq
 (D^\lm_{-j})_{mn}= (D^\lm_j)_{nm}
 \eeq
 we get
 \begin{eqnarray}\lb{41a}
 \sum_{s=0}^\infty (D^\lm_j)_{s s+j}U_{ms}U^*_{s+j n} &=&
 (D^\lm_{n-m})_{mn}t^\lm_{j n-m} \ \ \mathrm{for} \ j\geq 0 \\
 \sum_{s=0}^\infty (D^\lm_j)_{s-j s}U_{m s-j}U^*_{s n} &=&
 (D^\lm_{n-m})_{mn}t^\lm_{j n-m} \ \ \mathrm{for} \ j< 0,
 \end{eqnarray}
 \begin{eqnarray}\lb{42a}
 \sum_{s=0}^\infty t^\lm_{j s-m}(D^\lm_{s-m})_{ms} U_{sn} &=&
 (D^\lm_j)_{n-j n} U_{m n-j} \ \ \mathrm{for} \ n\geq j \\
 \sum_{s=0}^\infty t^\lm_{j s-m}(D^\lm_{s-m})_{ms} U_{sn} &=& 0 \ \
 \ \ \ \ \ \ \mathrm{for} \ n<j
 \end{eqnarray}
 and
 \beq\lb{43a}
 \sum_{s,l=0}^\infty t^\lm_{j l-s}(D^\lm_{l-s})_{s l}U^*_{ms}U_{l n}
 = (D^\lm_j)_{mn} \dl_{j n-m}.
 \eeq
 In the above formulas $U'$s and $t'$s are given in terms of the
 q-Kummer functions of the operator $\eta = B^*B$ (see (\ref{matru1}),
 (\ref{matru2}) and (\ref{matr}).

 In the coming section we give some simple examples.

 \vspace{1cm}
 \noindent
 {\bf 5. Examples}

 \vspace{2mm}
 \noindent
 {\bf A}. For  $n=m=0, \ j\geq 0$ (\ref{41a}) implies
 \beq
 \sum_{s=0}^\infty \fr{q^{\fr{s(1-s)}{2}}\eta^{2s} }{(s)_q!}
 \Phi^q(q^{-s}, q^{1+j}; \ q^{1+j+s}\lm^2) =
 \fr{(i\lm\eta)^{-j}(j)_q!}{\sqrt{e_q^{-\eta^2} e_q^{-q^j\eta^2}}}
 J^q_j(\lm \eta )
 \eeq
 which in $q\ar 1$ limit gives \cite{grad} (page 1038, Eq. 3 of
 8.975)
 \beq
 \sum_{s=0}^\infty \fr{\eta^{2s} }{s!}
 \Phi (-s, 1+j; \lm^2) =
 j! (\eta\lm)^{-j} e^{\eta^2}J_j(\lm \eta )
 \eeq
 {\bf B}. For $n=0$ and $k\equiv -j \geq m$ (\ref{42a}) implies
 \beq
 \sum_{s=0}^\infty
 q^{\fr{s(s+m-2k)}{2}}C_{sm}\eta^sJ_{s+k-m}^q(q^{\fr{s+k}{2}}\lm\eta)
 =
 \fr{q^{\fr{k(k-m)}{2}}\lm^k\eta^{k-m}}{(m)_q!(k-m)_q!}
 \Phi^q(q^{-m}, q^{1+k-m}; q^{1+k}\eta^2),
 \eeq
 where
 \beq
 C_{sm} = \fr{(-)^s \lm^{m-s}}{(s)_q!(m-s)_q!}\Phi_q (q^{-s},
 q^{1+m-s}; q^{1+m}\lm^2), \ \ \mathrm{for} \  m\geq s.
 \eeq
 For $s\geq m$ one has to replace $m, s$ with $s, m$ in the right
 hand side of the above expression. When $m=0$ we have
 \beq
 \sum_{s=0}^\infty
 q^{\fr{s^2}{2}+sk}\fr{(\lm\eta)^s}{(s)_q!}J^q_{s+k}(q^{\fr{s+k}{2}}\lm\eta)=
 q^{\fr{k^2}{2}}\fr{(\lm\eta)^k}{(k)_q!}
 \eeq
 which is the quantum analogue of a known formula \cite{grad} (page 974, Eq 1 of 8.515).

 \vspace{1mm}
 \noindent
 {\bf C}. For $j=\lm=0$ (\ref{43a}) implies the unitarity
 condition for the operator $U$:
 \beq
 \sum_{s=0}^\infty (U_{sm})^*U_{sn}= \dl_{nm},
 \eeq
 where we used
 \beq
 U^*_{ms} =(U_{sm})^*.
 \eeq
 For $n=m$ with $x=\eta^2$  we have
 \beq
 \sum_{s=0}^{n-1} \fr{q^{\fr{(n-s)(n-s+1)}{2}}(n)_q!x^{n-s}}{(s)_q!}(L_s^{q(s-n)}(x))^2 +
 \sum_{s=n}^{\infty} \fr{q^{\fr{(n-s)(n-s+1)}{2}}(s)_q!x^{s-n}}{(n)_q!}(L_n^{q(n-s)}(x))^2
 =e_{q^{-1}}^x
 \eeq
 In deriving the above examples one frequently uses the identities
 \beq
 B^kB^{*k}= q^{\fr{k(k+1)}{2}}\eta^{2k}, \ \ \ \
 B^{*k}B^k= q^{\fr{k(1-k)}{2}}\eta^{2k}.
 \eeq

 \end{document}